\documentclass[10pt,a4paper]{amsart}

\input xy
\xyoption{all}

\newtheorem{thm}{Theorem}[section]
\newtheorem{lemma}[thm]{Lemma}
\newtheorem{cor}[thm]{Corollary}
\newtheorem{pro}[thm]{Proposition}
\newtheorem{example}[thm]{Example}
\newtheorem{remark}[thm]{Remark}
\newtheorem{ddef}[thm]{Definition}

\def\ext{\mathop{\rm Ext}\nolimits}
\def\tot{\mathop{\rm Tot}\nolimits}
\def\tor{\mathop{\rm Tor}\nolimits}
\def\hom{\mathop{\rm Hom}\nolimits}
\def\spec{\mathop{\rm Spec}}
\def\proj{\mathop{\rm Proj}}
\def\codim{\mathop{\rm codim}}
\def\depth{\mathop{\rm depth}}
\def\reg{\mathop{\rm reg}\nolimits}

\def\b{{\beta}}

\def\fin{\mathop{\rm end}}
\def\ara{\mathop{\rm ara}}
\def\pd{\mathop{\rm pdim}}
\def\chara{\mathop{\rm char}}

\def\grade{\mathop{\rm grade}\nolimits}
\def\codim{\mathop{\rm codim}}
\def\indeg{\mathop{\rm indeg}}

\def\depth{\mathop{\rm depth}}

\def\supp{\mathop{\rm Supp}}
\def\sym{\mathop{\rm Sym}\nolimits}
\def\homgr{\mathop{\rm Homgr}\nolimits}

\def\C{{\mathcal C}}

\def\cc{{\C^{\bullet}_{\im}}}
\def\om{{\omega}}

\def\b{{\beta}}

\def\d{{\delta}}

\def\ip{{\mathfrak p}}

\def\im{{\mathfrak m}}
\def\cu{{\mathcal C}}

\def\cc{{\C^{\bullet}_{\im}}}

\def\ra{{\rightarrow}}
\def\lra{{\longrightarrow}}
\def\fini{{$\quad\quad\Box$}}

\newcommand{\bd}{\begin{ddef}}
\newcommand{\ed}{\end{ddef}}
\newcommand{\bt}{\begin{thm}}
\newcommand{\et}{\end{thm}}
\newcommand{\bl}{\begin{lemma}}
\newcommand{\el}{\end{lemma}}
\newcommand{\bco}{\begin{cor}}
\newcommand{\eco}{\end{cor}}
\newcommand{\bp}{\begin{pro}}
\newcommand{\ep}{\end{pro}}
\newcommand{\bex}{\begin{example}}
\newcommand{\eex}{\end{example}}
\newcommand{\brm}{\begin{remark}}
\newcommand{\erm}{\end{remark}}
\newcommand{\bconj}{\begin{conj}}
\newcommand{\econj}{\end{conj}}

\newcommand{\fm}{\frak{m}}

\newcommand{\fa}{\frak{a}}

\newcommand{\Ext}{\operatorname{Ext}}

\newcommand{\Tor}{\operatorname{Tor}}
\newcommand{\Hom}{\operatorname{Hom}}

\newcommand{\beqn}{\begin{eqnarray*}}
\newcommand{\eeqn}{\end{eqnarray*}}
\newcommand{\beq}{\begin{eqnarray}}
\newcommand{\eeq}{\end{eqnarray}}
\newcommand{\been}{\begin{enumerate}}
\newcommand{\eeen}{\end{enumerate}}

\begin{document}

\author[Chardin, Divaani-Aazar]{Marc Chardin and Kamran Divaani-Aazar}
\title[Generalized regularity]
{Generalized local cohomology and regularity of Ext modules}

\address{M. Chardin, Institut Mathematiques de Jussieu,
175 rue du Chevaleret,
F-75013, Paris, France.}
\email{chardin@math.jussieu.fr}

\address{K. Divaani-Aazar, Department of Mathematics, Az-Zahra
University, Vanak, Post Code 19834, Tehran, Iran.} \email{kdivaani@ipm.ir}



\begin{abstract} Let $R$ be a Noetherian  standard graded ring, and $M$ and
$N$ two finitely generated graded $R$-modules. We introduce
$\reg_R (M,N)$ by using the notion of generalized local cohomology instead 
of local cohomology, in the definition of regularity. We prove that $\reg_R (M,N)$
is finite in several cases. In the case that the base ring is a field, we show that
$$
\reg_R (M,N)=\reg (N)-\indeg (M).
$$
This formula, together with a graded version of duality for generalized local 
cohomology, gives a formula for the minimum of the initial degrees of  some Ext modules
(in the case $R$ is Cohen-Macaulay), of which the three usual definitions of regularity
are special cases. Bounds for regularity of certain Ext modules are obtained, using the 
same circle of ideas.
\end{abstract}

\maketitle

\section{Introduction}

Assume that $M$ is a graded module over a Noetherian  standard graded ring $R$. Then 
Castenluovo-Mumford regularity (or regularity for short) of $M$ is defined as
$$\reg (M):=\max_{i} \{\fin  (H_{R_+}^i(M))+i \}.$$

A generalization of local cohomology functors was given by Herzog in his Habilitationsschrift
\cite{H}. Let $\fa$ denote an ideal of a commutative Noetherian ring $R$. For each $i\geq 0$, 
 the functor $H^i_{\fa}(.,.)$ is defined by
$H^i_{\fa}(M,N)=\underset{n}{\varinjlim}\Ext_R^i(M/\fa^nM,N)$, for
all $R$-modules $M$ and $N$. This notion is a generalization of the usual 
local cohomology functor, that corresponds to the case $M=R$. This concept  has attained more notice
 in recent years, see e.g. \cite{HZ}. This notion inhierts many properties of the usual local 
cohomology. For example if $M$ and $N$ are  two finitely
generated $R$-modules such that $\pd M<\infty$, then, by \cite[2.5]{Y}, 
 $H^i_{\fa}(M,N)=0$ for all $i>\pd M+\ara (\fa)$, where
$\ara (\fa)$, the arithmetic rank of the ideal $\fa$, is the least
number of elements of $R$ required to generate an ideal which has
the same radical as $\fa$.

 Now, assume that $R=R_0 [R_1]$ is a Noetherian  standard graded ring and $M$ and
 $N$ are two finitely generated graded $R$-modules. Then, for any $i$,  $H_{R_+}^i(M,N)$ has a 
 natural graded structure, the $R_0$-modules $H_{R_+}^i(M,N)_n$  is zero for $n\gg 0$, and is finitely generated for all $n\in \Bbb Z$.

In view of the above mentioned properties of generalized local cohomology, it seems that
 there might be some connections between this notion and regularity. The main aim of the present paper is to investigate  such connections. Under the
above assumptions, we define the generalized regularity  of a pair $(M,N)$ of $R$-modules as 
$$\reg_{R} (M,N):=\sup_{i}\{\fin  (H_{R_+}^i(M,N))+i\}.$$ We prove that $\reg_{R} (M,N)$ is
finite in several cases.
If $R_0$ is a field and $M$ and $N$
are finitely generated, we show that  $$(*)\quad\reg_{R} (M,N)=\reg (N)-\indeg (M).$$

Notice that this shows in particular that $\reg_{R} (M,N)$ is independent of $R$ whenever
$R_0$ is a field. This is not the case for the modules $H^{i}_{R_+}(M,N)$, as for instance
$H^{i}_{R_+}(R_0,R_0)\simeq \ext^i_R(R_0,R_0)$ clearly depends on $R$.

In Sections 3 and 4, we will present several applications of this formula. In Section 3, we first 
establish a version of Grothendieck duality Theorem for graded generalized local cohomolgy.
Then by using $(*)$, we provide a sharp estimate for intial degrees of some $\Ext$ modules. As a 
corollary, we obtain a formula for the regularity which generalizes the three classical definitions
of Castelnuovo-Mumford regularity. Namely, we prove that if $R$ is a polynomial ring  over a field, 
then for any pair $M$, $N$ of non zero finitely generated graded $R$-modules
$$
\reg (M)-\indeg (N)=-\min\{ \indeg (\ext^{j}_{R}(M,N))+j\} .
$$

In Section 4, first by using formula $(*)$, an upper bound for the regluarity of the  $\Ext$ modules of 
certain pairs of modules is obtained. This result
contains a first result in this direction obtained by G. Caviglia in his thesis. 
Then we establish some upper bounds for the regularity of the  $\Ext$ modules 
of pairs of modules under an assumption on the dimension of some Tor modules. These results
are very much in the spirit of  the theorems 
obtained for the regularity of Tor modules, under essentially the same assumptions, by A. Conca 
and J. Herzog \cite{CH}, J. Sidman \cite{Si}, G. Caviglia \cite{Ca}, D. Eisenbud, C. Huneke and B. Ulrich \cite{EHU}, and by the first author \cite{Ch2}. \\

Throughout this paper, $R=R_0 [R_1]$ is a Noetherian standard graded ring, $\fm:=R_+$ and $M=\oplus 
_{n\in \Bbb Z} M_n$ and $N=\oplus_ {n\in \Bbb Z} N_n$  are two finitely generated graded 
$R$-modules. For a graded $R$-module $Q$, we denote by $\indeg (Q)$ (respectively $\fin (Q)$) the infimum 
(respectively the supremum) 
of the degrees of  non zero elements of $Q$ (with the convention $\indeg (0)=+\infty$ and $\fin (0)=-\infty$).

{\it We thank Juergen Herzog for his useful comments, and in particular for his suggestions 
regarding section 2.}

\section{Generalized regularity}

\bd
Let $a_{i}(M,N):=\fin (H^{i}_{\im}(M,N))$. Then,
$$
\reg_{R} (M,N):=\sup_{i}\{ a_{i}(M,N)+i\} .
$$ 
\ed

\bl\label{indegtor}
Let $f:R\lra S$ be a finite graded homomorphism of standard graded Noetherian rings. Assume 
that $R_0$ is local. 
Suppose $M$ is a finitely generated graded $S$-module and the following conditions hold :

i) there is a graded free $R$-resolution $F_{\bullet}^{S/R}$ of $S$ and an
integer $a$ such that $\indeg (F_{i}^{S/R})\geq i-a$, for all $i$,

ii) there is a graded free $S$-resolution $F_{\bullet}^{M/S}$ of $M$ and an
integer $b$ such that $\indeg (F_{i}^{M/S})\geq i-b$, for all $i$.

Then, a minimal graded free $R$-resolution $F_{\bullet}^{M/R}$ of $M$ satisfies $\indeg (F_{i}^{M/R})\geq i-a-b$, for all $i$.
\el

{\it Proof.} Let $k$ be the residue field of $R_0$. There is a spectral sequence
$$
E^{2}_{pq}=\tor_{p}^{S}(\tor_{q}^{R}(S,k),M)\ \Rightarrow\ \tor_{p+q}^{R}(k,M)\simeq F_{p+q}^{M/R}\otimes_{R}k.
$$
Assumption i) implies that $\tor_{q}^{R}(S,k)$ lives in degrees at least $q-a$. By
assumption ii), $E^{2}_{pq}$, which is a subquotient of $\tor_{q}^{R}(S,k)\otimes_{S}F_{p}^{M/S}$,
lives in degrees at least $(q-a)+(p-b)$. Hence $E^{\infty}_{pq}$ vanishes in degrees less than $(p+q)-a-b$.
The conclusion follows.\fini

\bco\label{regkk}
Let $R$ be a standard graded algebra over a regular local ring $(R_0,\im_0,k)$. 
The following are equivalent :

i) $\im_0$ contains an $R$-regular sequence of length $\dim R_0$,

ii) $R_\mu$ is a free $R_0$-module for any $\mu$.

If i) holds, then $\reg_{R}(k,k)=\dim R_0$. 
\eco

{\it Proof.} Notice that each $R_{\mu}$ is a finitely generated $R_0$-module of
finite projective dimension. A sequence of elements in $R_0$ is $R$-regular
if and only if it is $R_{\mu}$-regular for any $\mu$. Therefore condition i) is equivalent
to $\depth_{R_0}(R_{\mu})=\dim R_0$, which in turn is equivalent to ii). 

Assume that i) is satisfied. Let $S:=R/\im_0R$ and $f:R\lra S$ be the natural onto map. Set $d:=\dim R_0$.
Let $\pi :=(\pi_1,\ldots ,\pi_d)$ be an $R$-regular sequence generating $\im_0$.  The Koszul complex
$F_{\bullet}^{S/R}:=K_{\bullet}(\pi;R)$ is a graded free $R$-resolution of $S$,  $\indeg (F_{i}^{S/R})\geq i-d$, 
as $F_{i}^{S/R}=0$ for $i>d$. Let $F_{\bullet}^{k/S}$ and $F_{\bullet}^{k/R}$ denote minimal graded free $S$-resolution and $R$-resolution of $k$, respectively. As $S$ is a standard graded $k$-algebra, one has $\indeg (F_{i}^{k/S})\geq i$, for all $i$. By Lemma \ref{indegtor}, $\indeg (F_{i}^{k/R})\geq \max\{ 0,i-d\}$, which shows that
$H^{i}_{\im}(k,k)\simeq\ext^{i}_{R}(k,k)\simeq \hom_{R}(F_{i}^{k/R},k)$ vanishes in degrees above $d-i$, hence
$\reg_{R}(k,k)\leq d$.

On the other hand,  the degree 0 part of $F_{\bullet}^{k/R}$ is a minimal free $R_0$-resolution of $k$, so that
$\ext^{i}_{R}(k,k)_0\simeq \ext^{i}_{R_0}(k,k)\simeq \bigwedge^i k^d \not= 0$ for $i\leq d$. This proves
that $\reg_R(k,k)\geq d$.\fini

\bp\label{regkk2}
Let $f:S\lra R$ be a finite graded homomorphism of standard graded Noetherian rings. Assume that
$R_0$ and $S_0$ are local with the same residue field $k$, and that $f$ induces the identity on $k$. 
If $\indeg (\tor_i^S (R,k))>i$ for all $i>0$, then $\reg_R (k,k)\leq \reg_S (k,k)$.
\ep

{\it Proof.} We consider the change of ring spectral sequence,
$$
E^{2}_{pq}=\tor_{p}^{R}(\tor_{q}^{S}(R,k),k)\ \Rightarrow\ \tor_{p+q}^{S}(k,k).
$$
First notice that 
$$
\indeg (\tor_{p}^{R}(\tor_{q}^{S}(R,k),k))=\indeg (\tor_p^R(k,k))+\indeg (\tor_{q}^{S}(R,k))
$$
and therefore $\indeg (E^{2}_{pq})\geq \indeg (\tor_p^R(k,k))+q+1$ for $q>0$.

Recall that $\reg_S (k,k)=\max_i \{ i-\indeg (\tor_{i}^{S}(k,k))\}$, hence $\indeg (\tor_{i}^{S}(k,k)) \geq
i-\reg_S (k,k)$ for all $i\geq 0$. We may assume that $\reg_S (k,k)<\infty$. We use induction on $i$ to prove that 
$$
(*)\quad\quad \indeg (\tor_i^R(k,k))\geq i-\reg_S (k,k).
$$ 

For $i=0$, $\tor_0^R(k,k)=\tor_0^S(k,k)=k$. 

Now assume that $i>0$ and the inequality $(*)$
holds for any $j<i$. Let $\mu <i-\reg_S (k,k)$, then $\tor_{i}^{S}(k,k)_{\mu}=0$. As $R_0$ and $S_0$ are local with the same residue field $k$, and $f$ induces the identity on $k$, $\tor_{0}^{S}(R,k)=k\oplus N$,
for some finite $k$-module $N$. Therefore $(E^{2}_{i0})_{\mu}\simeq
\tor_{i}^{R}(k,k)_{\mu}\oplus N'$ for some finite $k$-module $N'$. Hence, it suffices to show that $(E^{2}_{i0})_{\mu}\simeq (E^{\infty}_{i0})_{\mu}$.
But for any $r\geq 2$, $E^{r+1}_{i0}\simeq \ker (E^{r}_{i0}\lra E^{r}_{i-r,r-1})$, and by induction,
$\indeg (E^{r}_{i-r,r-1})\geq \indeg (E^{2}_{i-r,r-1})\geq \indeg (\tor_{i-r}^R(k,k))+r\geq i-\reg_S (k,k)$. It follows that $(E^{r}_{i-r,r-1})_{\mu}=0$ for all  $r\geq 2$, which proves (*). Our claim follows since
$\reg_R (k,k)=\max_i \{ i+\fin (\ext^{i}_{R}(k,k))\}=\max_i \{ i-\indeg (\tor_{i}^{R}(k,k))\}$. \fini

\bt\label{regkksm}
Let $f:S\lra R$ be a surjective graded homomorphism of standard graded Noetherian rings. Assume that
$S_0$ is local with residue field $k$. Let $I:=\ker (f)$.

(1) If $\indeg (I)\geq \reg_S (k,k)+2$, then $\reg_R (k,k)\leq \reg_S (k,k)$.

(2) If $S$ is regular and $\indeg (I)\geq \dim (R_0)+1$, then $\reg_R (k,k)=\reg_S (k,k)=\dim R_0$.
\et

{\it Proof.} For (1), notice that for any $i\geq 0$, $\indeg (\tor_i^S(k,k))\geq i-\reg_S (k,k)$, and therefore 
$\indeg (\tor_i^S(I,k))\geq \indeg (I)+\indeg (\tor_i^S(k,k))\geq i+2$.
It implies that $\indeg (\tor_i^S(R,k))\geq i+1$ for any $i>0$ and the result follows from Proposition \ref{regkk2}.

 Now we prove (2). Our assumption on $\indeg (I)$ implies that $R_0=S_0$. Since $S$ is regular, we have $S\simeq R_0[X_1,\ldots ,X_n]$.  Set $X:=(X_1,\ldots ,X_n)$, $d:=\dim R_0$ and let $\pi :=(\pi_1,\ldots ,\pi_d)$ be an $R$-regular sequence generating $\im_0$.  The $i$-th homology group of the Koszul complex $K_{\bullet}:=K_{\bullet}(\pi ,X;I)$ is
$\tor_i^S(I,k)$. Notice that $(K_{i})_{\mu}=0$ for $\mu < \indeg (I)+\max\{ 0,i-d\}$, in particular, 
$\indeg (K_i)\geq \max\{ i,d\} +1$. Furthermore, if $i\geq d$, for $\mu = i-d+\indeg (I)$,
$$
H_i(K_{\bullet})_{\mu} =\ker (\bigwedge^{i-d}R_0^n\otimes_{R_0}I_{\indeg (I)}\lra (K_{i-1})_{i+1})
$$
and this kernel is zero since $\bigwedge^{i-d}R_0^n\otimes_{R_0}I_{\indeg (I)}\buildrel{\times \pi_1}\over{\lra}
\bigwedge^{i-d}R_0^n\otimes_{R_0}I_{\indeg (I)}\subset (K_{i-1})_{i+1}$ is already injective. Hence,
$\indeg (\tor_i^S(I,k))\geq i+2$ for any $i$. It follows that $\indeg (\tor_i^S(R,k))\geq i+1$ for any $i>0$.
Notice that $\reg_S (k,k)=\dim R_0$ by Corollary \ref{regkk} and that, for any standard graded Noetherain ring $T$ over a local base ring $(T_0 ,\im_0 ,k)$, $\reg _T (k,k)\geq \dim T_0$. 
Now (2) follows from Proposition \ref{regkk2}. 
\fini

\bl\label{chi}
Let $R$ be a standard graded Noetherian ring. Assume that
$R_0$ is a regular local ring of dimension $d>0$ with residue field $k$.
Then for any integer $\mu$,

(1) $\tor_i^R (k,k)_\mu =0$ for $i\geq (d+1)(\mu +1)$,

(2) $\sum_i (-1)^i \dim_k \tor_i^R (k,k)_\mu =0$.
\el

{\it Proof.} First notice that if $R$ is regular, then $R\simeq R_0[X_1,\ldots ,X_n]$, and (1) and 
(2) follows from the minimal free $R$-resolution of $k$ given by the Koszul complex which shows that
$\tor_i^R (k,k) \simeq \bigwedge^{i}(k^d\oplus  k[-1]^n)\simeq \oplus_\mu (\bigwedge^{i-\mu}k^d \otimes_k \bigwedge^{\mu}k^n)[-\mu ]$ which implies (2). It also shows that  $\indeg (\tor_i^R (k,k))=\max\{ 0,i-d\}$ for $0\leq i\leq d+n$ (and $+\infty$ else), which proves (1).

We induct on $\mu$. For $\mu =0$, (1) and (2) follows from the
isomorphisms $\tor_i^R (k,k)_0 \simeq \tor_i^{R_0} (k,k)\simeq \bigwedge^i k^d$.  

Let $S$ be a polynomial ring over $R_0$ with variables indexed by a set of minimal generators of the 
$R_0$-module $R_1$ and let $I:=\ker (S\ra R)$. We consider the change of ring spectral sequence,
$$
E^{2}_{pq}=\tor_{p}^{R}(\tor_{q}^{S}(R,k),k)\ \Rightarrow\ \tor_{p+q}^{S}(k,k).
$$
By the proof of Theorem \ref{regkksm}, $\indeg (\tor_{q}^{S}(R,k))\geq \max\{ 1,q-d+1\}$ for $q>0$.
Hence $\tor_{q}^{S}(R,k)\simeq \bigoplus_{j\geq \max\{ 1,q-d+1\}}k[-j]^{\b_{qj}}$, with $\b_{qj}:=
\dim_k \tor_{q}^{S}(R,k)_j$.

By induction hypothesis, $ (E^{2}_{pq})_\mu =0$ for $1\leq q\leq d$ and
$p\geq (d+1)\mu$ and for $q>d$ and $p\geq (d+1)(\mu -q+d)$. Hence
$ (E^{2}_{pq})_\mu =0$ for $p+q\geq (d+1)\mu +d$ and $q>0$. It follows that for
any $p>(d+1)\mu +d$, $ (E^{2}_{p0})_\mu \simeq \tor_{p}^{S}(k,k)_\mu =0$. This 
proves (1). 

For (2), the spectral sequence shows that 
$$
\sum_{p,q}(-1)^{p+q}\dim_k \tor_{p}^{R}(\tor_{q}^{S}(R,k),k)_\mu =\sum_{p,q}(-1)^{p+q}\dim_k (E^{\infty}_{pq})_\mu =
\sum_i (-1)^i \dim_k \tor_i^S (k,k)_\mu =0
$$
where the last equality follows from the regular case. It follows that
$$
\sum_{p}(-1)^{p}\dim_k \tor_{p}^{R}(k,k)_\mu +\sum_{\buildrel{q>0}\over{{\scriptscriptstyle{j\geq \max\{ 1,q-d+1\}}}}}(-1)^{q}
\b_{qj}\sum_p (-1)^p \dim_k \tor_{p}^{R}(k,k)_{\mu -j}=0.
$$
But $\sum_p (-1)^p \dim_k \tor_{p}^{R}(k,k)_{\mu -j}=0$ for $j>0$, by induction hypothesis. This proves (2).\fini

\brm
The bound in Lemma \ref{chi} (1) is sharp, as the following example shows. 
Let $k$ be a field, $A:=k[Y_1,\ldots ,Y_d,X_0,\ldots ,X_d]$ and $I$ be  the ideal of $2\times 2$ minors 
of $
\begin{pmatrix}
     X_0&X_1&\cdots&X_d\\
     0&Y_1&\cdots &Y_d\\
    \end{pmatrix}$. Consider
$R_0:=k[Y_1,\ldots ,Y_d]_{(Y_1,\ldots ,Y_d)}$,
$R:=R_0[X_0,\ldots ,X_d]/I'$, where $I'$ is the image of $I$ by the natural map from
$A$ to $S:=R_0[X_0,\ldots ,X_d]$.

By \cite[5.2.8]{A}, $B:=A/I$ is a Golod ring as it has minimal multiplicity. Notice that the minimal graded free
$R$-resolution of $k$ is the image of the minimal graded free $B$-resolution of $k$ via the natural map
from $B$ to $R$, as this map preserves exactness and sends the maximal graded ideal of $B$ to
the maximal graded ideal of $R$. Therefore, the graded Betti numbers of $k$ over $R$ are the same as the ones
of $k$ over $B$, graded by setting $\deg (X_i )=1$ and $\deg (Y_j )=0$. As $B$ is Golod, they are given
by the following formula for the Poincar\'e series (see \cite[(5.0.1)]{A}):
$$
P_{k}^{R}(t,u)={{P_k^A(t,u)}\over{1-t(P_{B}^A (t,u)-1)}}={{(1+t)^d (1+ut)^{d+1}}\over{1-t(\sum_{i=1}^{d}t^i{{d+1}\choose{i+1}}(\sum_{j=1}^{i}u^j))}},
$$
with $P_{M}^{T}(t,u):=\sum_{i,j}\dim_k \tor_i^T (M,k)_j t^i u^j\in {\bf Z}[u][[t]]$.

Indeed, $P_{k}^{A}(t,u)=(1+t)^d (1+ut)^{d+1}$ by the Koszul complex, and the minimal graded free $A$-resolution of $B$ is
given by the Eagon-Northcott complex $E_{\bullet}$, with $E_0=A$ and $E_i=(\wedge^{i+1}A^{d+1})[-1]\otimes_{A} \sym^{i-1}_{A}(A\oplus A[-1])$, which shows that $P_{B}^{A}(t,u)=1+\sum_{i=1}^{d}t^i{{d+1}\choose{i+1}}(\sum_{j=1}^{i}u^j)$. In particular,
$P_{k}^{R}(t,u)$ has the monomials $t^d(t^{d+1}u)^j$ in its expension, which shows that $\dim_k \tor_{d+(d+1)j}^{R}(k,k)_j=1$, and the sharpness of the bound.
\erm

\brm
i) The conditions i) and ii) in Corollary \ref{regkk} are satisfied when $R$ is regular ({\it i.e.} when
$R$ is a polynomial ring over a regular local ring). 

ii) Notice that $\reg_R (k,k)$ is infinite whenever $R_0$ is not regular. The following example illustrates
that even when $R_0$ is a discrete valuation ring,  $\reg_R (k,k)$ can be infinite. 

iii) The following example illustrates that generalized local cohomology, unlike local 
cohomology, is very dependent on R. It shows this dependence for generalized regularity
as well.
\erm

\bex\label{piX}
Let  $(R_0,\pi )$ be a discrete valuation ring, $k:=R_0/(\pi )$ and $R:=R_0[X]/(\pi X)$. Then $R$ is a 
standard graded algebra over $R_0$, which is local  regular  of dimension 1. Nevertheless
$\reg_R (k,k)=\infty$. Indeed the minimal graded free $R$-resolution $F_{\bullet}$ of $k$ is
$$
\xymatrix{
\cdots\ar^{\psi}[r]&F_{2i+1}\ar^{\phi}[r]&F_{2i}\ar^{\psi}[r]&F_{2i-1}\ar^{\phi}[r]&\cdots\ar^{\psi}[r]&F_{1}
\ar[r]&R\\}
$$
with $F_{2i}:=R[-i]^2$ and $F_{2i-1}:=R[-i+1]\oplus R[-i]$ for $i>0$, 
$\psi :=
\begin{pmatrix}
     x&0\\
     0&\pi 
     \\
    \end{pmatrix}$ and
$\phi :=\begin{pmatrix}
     \pi &0 \\
     0&x\\
    \end{pmatrix}$.
This resolution shows that for $i>0$, $\ext_R^{2i-1}(k,k)$ sits in degrees
$-i$ and $-i+1$ and $\ext_R^{2i}(k,k)$ sits in degree $-i$. 

Now, $H^j_\im (k,k)\simeq \ext_R^{j}(k,k)$, and therefore $a_{j}(k,k)=-[j/2]$ so that
$a_{j}(k,k)+j=[(j+1)/2]$ is unbounded.
\eex

\bt
Let $R$ be a standard graded Noetherian ring. Assume that
$R_0$ is a discrete valuation ring with residue field $k$.
Then the following are equivalent :

(i) $R_1$ is a free $R_0$-module,

(ii) $\reg_R (k,k)=1$,

(iii) $\reg_R (k,k)<\infty$,


(iv) $\reg_R (k,M)=\reg (M)+1$, for any $M$.

Moreover if one of these conditions does not hold, then $\indeg (\tor_i^R(k,k))=\lfloor i/2\rfloor$ for all $i\geq 0$.
\et

{\it Proof.} First $(i)\Rightarrow (ii)$ by Theorem \ref{regkksm} (2) and clearly $(iv)\Rightarrow (ii)\Rightarrow (iii)$. For $(iii)\Rightarrow (i)$, we will show that $\tor_{2i+1}^R(k,k)_i\not= 0$ for all $i\geq 0$ when $R_1$
is not a free $R_0$-module. 
Let $S=R_0 [X_1,\ldots ,X_n ]$ be a polynomial ring over $R_0$ with variables indexed by a set of minimal generators of the 
$R_0$-module $R_1$ and let $I:=\ker (S\ra R)$. Choose $f=\pi^m l\in I_1\setminus \pi I_1$, where $m\geq 1$ and $l\not\in (\pi )$. Set $S':=S/(f)$ and $I':=IS'$. We consider the change of ring spectral sequence,
$$
E^{2}_{pq}=\tor_{p}^{R}(\tor_{q}^{S'}(R,k),k)\ \Rightarrow\ \tor_{p+q}^{S'}(k,k).
$$

A minimal graded free $S'$-resolution $F_\bullet$ of $S'/(l,\pi )S'$ is given by
$$
\xymatrix{
\cdots\ar^{d_{2i+1}}[r]&F_{2i}\ar^{d_{2i}}[r]&F_{2i-1}\ar^{d_{2i-1}}[r]&\cdots\ar^{d_{3}}[r]&F_{2}\ar^(.4){d_{2}}[r]&S'\oplus S'[-1]
\ar^(.6){(\pi\ l)}[r]&S'\\}
$$
with $F_{2i}:=S'[-i]^2$ and $F_{2i+1}:=S'[-i-1]\oplus S'[-i]$ for $i>0$, 
$d_2 :=\begin{pmatrix}
     l&0\\
     -\pi&\pi^m \\
    \end{pmatrix}$,
 $d_3 :=\begin{pmatrix}
     0&\pi^m \\
     l&\pi \\
    \end{pmatrix}$,
and for $i\geq 2$,
$d_{2i} :=\begin{pmatrix}
     \pi &0\\
     -l&l\pi^{m-1} \\
    \end{pmatrix}$,
$d_{2i+1} :=\begin{pmatrix}
     l\pi^{m-1}&0\\
    l &\pi\\
    \end{pmatrix}$.

As $l=\sum_i c_i X_i \not\in \pi I_1$, there exits $i$ with $c_i\not\in (\pi )$, and a minimal free $S'$-resolution of $k$ is given by $G_\bullet :=F_{\bullet}\otimes_{S'}K_{\bullet}(X_1,\ldots ,\widehat{X_i},\ldots  ,X_n;S')$ (recall that $S'/(l,\pi )S'\simeq k[X_1,\ldots ,\widehat{X_i},\ldots  ,X_n]$). One has $G_i=\oplus_{j=0}^{i}\left( F_{i-j}\otimes_{S'}S'[-j]^{{{n-1}\choose{j}}}\right) $, in particular
$\indeg (G_i)=\indeg (F_i)$ for all $i$.

It follows that, for $i\geq 0$,  $\indeg (\tor_{2i+1}^{S'}(R,k))\geq\indeg (\tor_{2i}^{S'}(I',k))\geq i+1$.
 Furthermore, $\tor_{2i-1}^{S'}(I',k)_i\simeq \ker (I'_1\buildrel{\times \pi}\over{\lra} I'_1)=0$,   because 
$(lS_0\cap I_1)/fS_0=0$.
It follows that $\indeg (\tor_{2i}^{S'}(R,k))=\indeg (\tor_{2i-1}^{S'}(I',k))\geq i+1$ and  
therefore $\indeg (\tor_{q}^{S'}(R,k))\geq \lfloor q/2\rfloor +1$ for $q>0$.

For $q>0$ we get from Lemma \ref{chi} and the above estimates that 
$$
\indeg (E^{2}_{pq})=\indeg (\tor_{p}^{R}(k,k))+\indeg (\tor_{q}^{S'}(R,k))\geq \lfloor p/2\rfloor +  \lfloor q/2\rfloor +1 = i+1
$$
for $p+q=2i+1$. Hence, for any $i\geq 0$, $(E^{\infty}_{2i+1,0})_i\simeq\tor_{2i+1}^{S'}(k,k)_i\not= 0$  is a submodule of $(E^{2}_{2i+1,0})_i\simeq \tor_{2i+1}^R(k,k)_i$, which is therefore not zero. 
This proves $(iii)\Rightarrow (i)$ and also shows that if $(i)$ is not satisfied, then $\indeg (\tor_{2i+1}^{R}(k,k))\leq i$.
In this case, $\indeg (\tor_{2i+1}^{R}(k,k))=i$  by Lemma \ref{chi}. Furthermore $\indeg (\tor_{2i}^{R}(k,k))=i$ as 
$\indeg (\tor_{2i}^{R}(k,k))\leq  \indeg (\tor_{2i+1}^{R}(k,k))=i$ on one side and $\indeg (\tor_{2i}^{R}(k,k))\geq i$ by
Lemma \ref{chi} on the other side.

Now $(ii)\Rightarrow (iv)$ by Theorem \ref{greg} 3). 
\fini

As the following example shows, the finiteness of generalized regularity is not 
symmetric in $M$ and $N$.

\bex
Let  $(R_0,\pi )$ be a discrete valuation ring, $k:=R_0/(\pi )$ and $R:=R_0[X]/(\pi X, X^2)$. 
Then $R$ is a 
standard graded algebra over $R_0$, which is local  regular  of dimension 1, $\reg_R (R,k)=0$ and
$\reg_R (k,R)=\infty$. Indeed the minimal graded free $R$-resolution $F_{\bullet}$ of $k$ is given by
$$
\xymatrix{
\cdots\ar^(.3){d_3\oplus d_2}[r]&F_3=(F_{2}\oplus F_{1})[-1] \ar^(.5){d_2\oplus d_1}[r]&F_2=(F_{1}\oplus F_{0})[-1] \ar^(.7){d_2}[r]&F_{1}
\ar^{d_1}[r]&F_0\\}
$$
with  $F_0:=R$, $F_1:=R\oplus R[-1]$, $d_1$ is given by the matrix $\begin{pmatrix}x&\pi \end{pmatrix}$ and
$d_2$ is given by the matrix $\begin{pmatrix}x&0&\pi \\ 0&x&0\end{pmatrix}$. The first three steps of the resolution
were computed using the software Macaulay 2 \cite{GS}, and the rest of the resolution follows
by induction, due to the decomposition of the map from $F_3$ to $F_2$. The decomposition of the
modules and maps in the resolution also shows that for any module $M$ and any $i\geq 3$,
$$
\ext^i_{R}(k,M)\simeq \ext^{i-1}_{R}(k,M)[1]\oplus \ext^{i-2}_{R}(k,M)[1].
$$
Using  Macaulay 2, it is easy to check that $\fin (\ext^i_{R}(k,R))=0$ 
for $i=0,1,2$. It follows that $\fin (\ext^i_{R}(k,R))=-[(i-1)/2]$ for any $i>0$, which proves that 
$\reg_R (k,R)=\infty$.
    \eex


If $F_\bullet$ and
$G^\bullet$ are two complexes of $R$-modules, $\homgr_{R}(F_{\bullet}, G^{\bullet})$ is the cohomological complex
with modules $C^i=\prod_{p+q=i}\hom_R(F_p,G^q)$. If either $F_{\bullet}$ or $G^{\bullet}$ is bounded,
 or both $F_{\bullet}$  and $G^{\bullet}$ are bounded below, then
$\homgr_{R}(F_{\bullet}, G^{\bullet})$ is the totalisation of the double complex with $C^{p,q}=\hom_R(F_p,G^q)$.

We denote by $\cc (\hbox{---})$ the \v Cech complex on $\hbox{---}$.

 \bl\label{GCech}
 Let $F_{\bullet}^{M}$ be a graded free $R$-resolution of $M$ and $I^{\bullet}_{N}$ be a graded
 injective $R$-resolution of $N$. Consider $C^{\bullet\bullet}$, the double complex with $C^{p,q}=\hom_{R}(F_{q}^{M},\cu_\im^{p}(N))$,
$D^{\bullet\bullet}$  the double complex with $D^{p,q}=\hom_{R}(F_{q}^{M},\Gamma_{\im}I^{p}_{N})$, $T_{C}^{\bullet}:=\tot (C^{\bullet\bullet})$
and $T_{D}^{\bullet}:=\tot (D^{\bullet\bullet})$. Then,
$$
 H^{i}_{\im}(M,N)\simeq H^{i}(T_{C}^{\bullet})\simeq H^i (\homgr_R (F_{\bullet}^{M},\cc (N))\simeq H^{i}(T_{D}^{\bullet})\simeq H^i (\homgr_R (F_{\bullet}^{M},\Gamma_{\im}I^{\bullet}_{N})).
 $$
 \el
 
 {\it Proof.} All these modules are the total homology of the triple complex $T^{\bullet \bullet \bullet }$ with $T^{i,j,k}= \hom_{R}(F_{i}^{M},\cu_\im^{j}(I^{k}_{N}))$. On the other hand, $H^{i}(T_{D}^{\bullet})\simeq H^{i}(\hom_{R}(M,\Gamma_{\im}I^{\bullet}_{N}))\simeq  H^{i}_{\im}(M,N)$. \fini
 
 \bl\label{ssHE} Let $F_{\bullet}$ be a free $R$-resolution of $M$. There exist two convergent spectral sequences, 
$$
{'E}_{1}^{pq}=H^{p}_{\im}(N)\otimes_{R} \hom_{R}(F_{q},R)\ \Rightarrow \ H^{p+q}_{\im}(M,N)
$$
and
$$
{''E}_{2}^{pq}=H^{p}_{\im}(\ext^{q}_{R}(M,N))\ \Rightarrow \ H^{p+q}_{\im}(M,N).
$$

In particular, $a_{i}(M,N)$ is finite for any $i$, and $\reg_{R}(M,N)$ is finite if $\ext^i_R (M,N)=0$ for $i\gg 0$.
\el

{\it Proof.} Consider the double complex with 
$$
C^{pq}:=\hom_{R}(F_{q},\C^{p}_\im (N))\simeq \hom_{R}(F_{q},R)\otimes_{R} \C^{p}_\im (N)
\simeq \C^{p}_\im (\hom_{R}(F_{q},N)).
$$
This complex has total homology  isomorphic to $H^{\bullet}_{\im}(M,N)$ by Lemma \ref{GCech}, and first or second terms of the spectral sequences associated to row and column filtrations are as stated.\fini

The following easy observation will be useful for the proof of the last theorem of this section.

\bl\label{endhom} Let  $P$ and $Q$ be two graded $R$-modules. 
If $\indeg (P)$  and $\fin (Q)$ are finite, then $$\fin  ({^*\Hom}_R(P,Q))\leq \fin  (Q)-\indeg (P).$$ 
Equality holds 
if $P$ is free of finite rank.
\el

In the sequel, for a finitely generated graded $R$-module $P$, we set $a_i(P):=\fin (H^i_\im (P))$.

\begin{thm}\label{greg} 
Assume that $(R_0,\im_0 ,k)$ is local.\\
1) For any $i$, $a_{i}(M,N)\leq \reg (N)-\indeg (M)$ and if $R_0$ is a field then $a_{i}(M,N)\leq \reg (N)-\indeg (M)-i$.\\
2) $0\leq \reg_{R}(M,k)+\indeg (M)\leq \reg_{R} (k,k)$  for any $M\not= 0$.\\
3) If $\reg_{R} (M,k)$ is finite, then 
$$
\reg_{R} (M,N)=\reg (N)+\reg_{R} (M,k).
$$
4) If $\reg_{R} (M,k)$ is finite, $\reg_{R} (M,N)=a_p (M,N)+p$ for 
$$
p:=\min \{ i\ \vert\ \reg_{R} (M,k)=a_i (M,k)+i\} +  \max\{ j\ \vert\ \reg (N)=a_j (N)+j\} .
$$
5) If $R_{0}$ is a field, then $\reg_{R} (M,N)= \reg (N)-\indeg (M)=a_p (M,N)+p$ for any $p$ such that 
$\reg (N)=a_p (N)+p$.\\
\end{thm}

{\it Proof.} Let $F_{\bullet}$ be a minimal graded free $R$-resolution of $M$ and set $C^{pq}:=\C^{p}_{\im}(\hom_{R}(F_{q},N))$ as in Lemma \ref{ssHE}.
Set $e:=\indeg (M)$ and $F^{\bullet}:=\hom_{R}(F_{\bullet},R)$. 
Notice that $C^{pq}\simeq F^{q}\otimes_{R}\C^{p}_{\im}(N)$. As the second spectral sequence of Lemma \ref{ssHE} provides an isomorphism $H^q_{\im}(M,N)\simeq \ext^q_R(M,N)$ whenever
$M\otimes_R N$ is supported in $\im$, one has
$$
a_{i}(M,k)=\fin (\ext^{i}_{R}(M,k))=\fin (F^i\otimes_{R}k)=-\indeg (F_i).
$$

One has $('E_{1}^{pq})_{\mu}=\hom_{R}(F^{q},H^{p}_{\im}(N))_{\mu}=0$
for $\mu >a_{p}(N)-e$, by Lemma \ref{endhom}. 
If further $R_0$ is a field,  $\indeg (F_{j+1})>\indeg (F_{j})$ for any $j$ and Lemma \ref{endhom} implies that $('E_{1}^{pq})_{\mu}=0$ for $\mu >a_{p}(N)-e-q$. 
 This proves 1). 

If $\reg_{R}(k,k)$ is finite, a minimal free $R$-resolution $F_{\bullet}^k$ of $k$ satisfies
$\indeg (F_i^k)\geq i-c$, for  $c:=\reg_R (k,k)$. It follows that $\indeg (F_i)=\indeg (\tor_i^R(M,k))\geq \indeg (M)+i-c$,
as $\tor_i^R(M,k)\simeq H_i(M\otimes_R F_{\bullet}^k)$. This proves 2).

By \cite[5.1 (ii)]{Ch2},  if $\reg_{R} (M,k)=\sup_{i}\{ \fin (F^i\otimes_{R}k)+i\}<\infty$,
then setting $i_{0}=\min \{ i\ \vert\ \reg_{R} (M,k)=a_i (M,k)+i\}$ and $j_{0}:=\max\{ j\ \vert\ \reg (N)=a_j (N)+j\}$,
one has :
$$
\reg_{R}(M,N)=a_{i_0+j_0}(M,N)+(i_0+j_0)=\reg_{R}(M,k)+\reg (N).
$$
This proves 3) and 4).

For 5), first notice that 1) implies that $\reg_R (M,N)\leq \reg (N)-\indeg (M)$. 
Set $\mu :=a_{p}(N)-e=\reg (N)-p-e$. 
One has $'E_{r}^{p+r-1,-r}=0$ for any $r\geq 1$,  so that $'E_{r+1}^{p0}$ is isomorphic to the kernel of the map
$$
\xymatrix{
'E_{r}^{p0}\ar^(.45){{'d}_{r}^{p0}}[rr]
&&'E_{r}^{p-r+1,r}\\}
$$
But $('E_{1}^{p-r+1,r})_{\mu}\simeq \hom_{R}(F_{r},H^{p-r+1}_{\im}(N))_{\mu}$  
is zero for $r\geq 1$ because
$$
\fin (\hom_{R}(F^{r},H^{p-r+1}_{\im}(N)))
\leq (\reg (N)-p+r-1)+(-e-r)
\leq \mu -1.
$$
Therefore $('E_{\infty}^{p0})_{\mu}\simeq (H^{p}_{\im}(N)\otimes_{R}F^{0})_{\mu}
\simeq H^{p}_{\im}(N)_{a_{p}(N)}^{\dim_{k}(M_{e})}\not= 0$. It follows that $H^{p}_{\im}(M,N)_{\mu}\not= 0$, 
which proves that $\reg_{R} (M,N)\geq \mu +p=\reg (N)-\indeg (M)$.\fini

\brm \label{vanfpd}
i) If $R_0$ is regular, not necessarily local, and either $R_\mu$ is a free $R_0$-module 
for any $\mu$ or $R$ is the quotient of a polynomial ring over $R_0$ by an ideal of initial
degree at least $\dim R_0+1$, then 
$$
\reg_R (M,N)\leq \reg (N)-\indeg (M)+\dim R_0.
$$
This follows from the fact that $\reg_R (M,N)=\max_{\ip_0\in \spec (R_0)}\{ \reg_{R_{\ip_0}} (M_{\ip_0},N_{\ip_0})\}$, and that for any $\ip_0\in \spec (R_0)$, $\reg (N_{\ip_0})\leq \reg (N)$ and 
$\indeg (M_{\ip_0})\geq \indeg (M)$.

ii) If $M$ has finite projective dimension, or $N$ has finite injective dimension, then $\ext^{q}_{R}(M,N)$ has
dimension at most $\dim (R)-q$ for any $q$. Hence,
the spectral sequence $''E_2^{pq}=H^{p}_{\im}(\ext^{q}_{R}(M,N))\Rightarrow H^{p+q}_{\im}(M,N)$ 
shows that  $H^{i}_{\im}(M,N)=0$ for $i>\dim (R)$.
\erm

 \section{Duality for generalized local cohomology}
 
We establish in this section a duality for graded generalized local cohomology in our context, which is a 
graded version of the original result of Suzuki \cite{Su} (see also \cite{HZ} for further results).
 
 \bt\label{duality}
 Let $R$ be a standard graded Cohen-Macaulay ring of dimension $d$. Assume that
 $R_{0}$ is a field and let $\om_{R}$ be the canonical module of $R$, and $M$ and $N$ be finitely generated 
 graded $R$-modules. 
 
 1) $H^{i}_{\im}(M,R)\simeq {^{*}\hom}_{R}(\tor_{i-d}^{R}(M,\om_{R}),R_0)$.

 2) If $M$ has finite projective dimension, or if $N$ has finite projective dimension and
$\tor_{i}^{R}(M,\om_{R})=0$ for $i>0$,  then
 $$
 H^{i}_{\im}(M,N)\simeq {^*\hom}_{R}(\ext^{d-i}_{R}(N,M\otimes_{R}\om_{R}),R_{0}).
 $$
 
 \et

{\it Proof.} 
Let $F_{\bullet}^{M}$ be a minimal graded free $R$-resolution of $M$ and $F_{\bullet}^{N}$  be a minimal graded free $R$-resolution of $N$.
By Lemma \ref{GCech} and \cite[\S 6, Th\'eor\`eme 1, b]{B},
$$
H^{i}_{\im}(M,N)\simeq H^{i}(\homgr_{R}(F_{\bullet}^{M},\cc (N)))\simeq 
H^{i}(\homgr_{R}(F_{\bullet}^{M},\cc (F_{\bullet}^{N})).
$$

As $R$ is Cohen-Macaulay, $H^{i}_{\im}(R)=0$ for $i\not= d$, and therefore
$$
H^{i}_{\im}(M,N)\simeq H^{i-d}(\homgr_{R}(F_{\bullet}^{M}, F_{\bullet}^{N}\otimes_{R}H^{d}_{\im}(R))).
$$
As $F_{\bullet}^{M}$ or $F_{\bullet}^{N}$ is bounded 
$$
\begin{array}{rcl}
H^{i}_{\im}(M,N)&\simeq&H_{d-i}(\hom_{R}(\homgr_{R}(F_{\bullet}^{N},F_{\bullet}^{M}),H^{d}_{\im}(R)))\\
&\simeq&H_{d-i}({^{*}\hom}_{R}(\homgr_{R}(F_{\bullet}^{N},F_{\bullet}^{M}\otimes_{R}\om_{R}),R_0))\\
&\simeq&{^{*}\hom}_{R}(H^{d-i}(\homgr_{R}(F_{\bullet}^{N},F_{\bullet}^{M}\otimes_{R}\om_{R})),R_0).\\
\end{array}
$$
This shows 1) and, by Remark \ref{vanfpd} ii), it implies that $\tor_{i}^{R}(M,\om_{R})=0$ for $i>0$ if $F_{\bullet}^{M}$ is bounded.  Hence $F_{\bullet}^{M}\otimes_{R}\om_{R}$ is acyclic, and 2) follows.
\fini

\begin{thm}
Let $R$ be a standard graded Gorenstein ring. Assume that
 $R_{0}$ is a field, $M$ and $N$ are finitely generated graded $R$-modules and $M$ or $N$ 
 has finite projective dimension, then
$$
\min_j \{ \indeg (\ext^{j}_{R}(N,M))+j\} =\reg (R)-\reg (N)+\indeg (M).
$$
\end{thm}

{\it Proof.} Let  $d:=\dim R$ and $a$ be the $a$-invariant  of $R$. By Theorem \ref{duality},
$$
H^{i}_{\im}(M,N)\simeq (\ext^{d-i}_{R}(N,M)[a])^{*}.
$$
In particular,
$$
-a_{i}(M,N)=\indeg (\ext^{d-i}_{R}(N,M)[a])=\indeg (\ext^{d-i}_{R}(N,M))-a.
$$
Therefore
$$
\begin{array}{rcl}
\reg_{R} (M,N)&=&\max_{i}\{ -\indeg (\ext^{d-i}_{R}(N,M))+a+i\}\\
&=&a+d-\min_{j}\{ \indeg (\ext^{j}_{R}(N,M))+j\}.
\end{array}
$$
The conclusion follows from Theorem \ref{greg} 5), since $a+d=\reg (R)$.\fini

\begin{cor}
If $R$ is a polynomial ring  over a field, then for any pair of non zero finitely generated 
graded $R$-modules $M$ and $N$,
$$
\reg (M)-\indeg (N)=-\min\{ \indeg (\ext^{j}_{R}(M,N))+j\} .
$$
\end{cor}

Applying this corollary to $N:=R_0$ we recover the definition of regularity in terms of graded 
Betti numbers, applying it with $N:=R$ we recover the definition of regularity as in Eisenbud's
book [Ei, 20.16], which is equivalent to the one in terms of local cohomology by graded local duality.

\section{Regularity of Ext modules}

The following result generalizes a result of Caviglia in his thesis (see \cite[3.10]{Ca}).

\bt Let $R$ be a Noetherian standard graded ring of dimension $n$
such that $R_0$ is a field. Let $M$ and $N$ be two finitely generated graded 
$R$-modules and $i_0$ an integer. Assume that $\Ext_R^i(M,N)$ is zero for $i<i_0$,
$\Ext_R^{i_0}(M,N)$ has dimension at most $n-i_0$ and $\Ext_R^i(M,N)$ is either
a Cohen-Macaulay $R$-module of dimension $n-i$ or zero for all $i>i_0$. Then,\\
i) $H_{\fm}^{i-i_0}(\ext^{i_0}_R(M,N))\simeq H_{\fm}^{i}(M,N)$ for all $i\not= n$. \\
ii) $\max_{i} \{\reg (\Ext_R^i(M,N))+i \}=\reg (N)-\indeg (M)$. \\
iii) $\reg (\ext^{i_0}_R(M,N))+i_0=\reg (N)-\indeg (M)$ if $\reg (N)=a_{p}(N)+p$
for some $p<n$.
\et

Recall that $\Ext_R^{i_0}(M,N)$ has dimension at most $n-i_0$ if either $M$  has finite projective dimension or $N$ has finite injective dimension (see Remark \ref{vanfpd} ii)).

{\it Proof.} Lemma \ref{ssHE} provides a spectral sequence
$$
{''E}_{2}^{pq}=H^{p}_{\im}(\ext^{q}_{R}(M,N))\ \Rightarrow\ H^{p+q}_{\im}(M,N)
$$
and our hypothesis implies that ${''E}_{2}^{pq}=0$ for $p\not= n-q$ and $q\not= i_0$, and for $p+q>n$. This in turn shows 
that ${''E}_{2}^{pq}\simeq {''E}_{\infty}^{pq}$ for any $p$ and $q$. It proves i) and provides a 
filtration $0=F_{-1}\subseteq F_{0}\subseteq \cdots \subseteq F_{n-i_0}=H^{n}_{\im}(M,N)$ with
$F_i/F_{i-1}\simeq H^{i}_{\im}(\ext^{n-i}_{R}(M,N))$. It follows that
$$
a_{n}(M,N)=\max_{i}\{ a_{i}(\ext^{n-i}_{R}(M,N))\} .
$$
 Notice that $a_{i}(\ext^{n-i}_{R}(M,N))=\reg (\ext^{n-i}_{R}(M,N))-i$ for
$i\not= n-i_0$. Hence, by i), 
$$
\reg_{R} (M,N)=\max\{ \max_{i}\{ a_{i-i_0}(\ext^{i_0}_{R}(M,N))+i\} ,\max_{i\not= n-i_0}\{ \reg (\ext^{n-i}_{R}(M,N))+n-i\}\} ,
$$
which proves ii) by Theorem \ref{greg} 5). 




If  $\reg (N)=a_{p}(N)+p$ for some $p<n$, then $\reg_{R} (M,N)=a_p (M,N)+p$ by Theorem \ref{greg} 5). Hence by i), $\reg_{R} (M,N)=a_{p-i_0}(\ext^{i_0}_{R}(M,N))+p\leq 
\reg (\ext^{i_0}_{R}(M,N))+i_0$, which proves iii) in view of ii).\fini

\begin{lemma}\label{ssTE} If $M$ or $N$ has finite projective dimension, there exists a spectral sequence :
$$
\tor_{i}^{R}(\ext^{j}_{R}(M,R),N)\Rightarrow \ext^{j-i}_{R}(M,N).
$$
\end{lemma}

{\it Proof.} Let $F_{\bullet}^{M}$ and $F_{\bullet}^{N}$ be minimal graded free $R$-resolution of 
$M$ and $N$, respectively. The double complex with terms $\hom_{R}(F_{p}^{M},F_{q}^{N})\simeq \hom_{R}(F_{p}^{M},R)\otimes_{R}F_{q}^{N}$ has total homology isomorphic to  
$\ext^{\bullet}_{R}(M,N)$, because $M$ or $N$ has finite projective dimension. It gives rise
to a spectral sequence with $E_2$ term as in the Lemma.\fini

\begin{cor}\label{mapExt} There exists a natural  map 
$$
\ext^{j}_{R}(M,R)\otimes_{R}N\lra \ext^{j}_{R}(M,N),
$$
whose kernel is supported on  $\bigcup_{i\geq 1}\supp (\tor_{i+1}^{R}(\ext^{j+i}_{R}(M,R),N))$
and whose cokernel is supported on $\bigcup_{i\geq 1}\supp (\tor_{i}^{R}(\ext^{j+i}_{R}(M,R),N))$.
\end{cor}

{\it Proof.} The spectral sequence of Lemma \ref{ssTE} provides a natural onto map from $\ext^{j}_{R}(M,R)\otimes_{R}N$ to ${^{h}_{\infty}E}_{0}^{j}$ whose kernel is supported on  $\bigcup_{i\geq 1}\supp (\tor_{i+1}^{R}(\ext^{j+i}_{R}(M,R),N))$ and a natural into map from ${^{h}_{\infty}E}_{0}^{j}$ to $\ext^{j}_{R}(M,N)$ whose 
cokernel is supported on $\bigcup_{i\geq 1}\supp ({^{h}_{\infty}E}_{i}^{j+i})$. The conclusion follows as 
${^{h}_{\infty}E}_{i}^{j+i}$ is a subquotient of  ${^{h}_{2}E}_{i}^{j+i}\simeq \tor_{i}^{R}(\ext^{j+i}_{R}(M,R),N)$.\fini\medskip

We  assume from this point on that $R=R_0 [X_1,\ldots ,X_n]$ is a polynomial ring over a field $R_0$.

\begin{cor}\label{reg2E} Let $c$ and $e$ be two integers. Assume that $M_{\ip}$ is Cohen-Macaulay of 
codimension $c$ for all $\ip \in \supp (N)$ such that $\dim R/\ip >e$, then
$$
H^{i}_{\im}(\ext^{c}_{R}(M,R)\otimes_{R}N)\simeq H^{i}_{\im}(\ext^{c}_{R}(M,N)),\ \forall i\geq e+2,
$$
and there is an onto map
$$
H^{e+1}_{\im}(\ext^{c}_{R}(M,R)\otimes_{R}N)\lra H^{e+1}_{\im}(\ext^{c}_{R}(M,N)).
$$
\end{cor}
 
 {\it Proof.} The assumption implies that $\ext^{j}_{R}(M,R)_{\ip}=0$  if $\ip \in \supp (N)$, $\dim R/\ip >e$ and $j\not= c$ because $H^{\dim R_\ip -j}_{\ip}(M_\ip )=0$ in these cases. This in turn shows that $\tor_{i}^{R}(\ext^{c+i}_{R}(M,R),N)$ and $\tor_{i+1}^{R}(\ext^{c+i}_{R}(M,R),N)$ are supported in dimension at most $e$ for any $i\geq 1$. By Corollary \ref{mapExt}, there is an exact sequence
 $$
0\lra K\lra \ext^{c}_{R}(M,R)\otimes_{R}N\lra \ext^{c}_{R}(M,N)\lra C\lra 0,
$$
with $\dim (K)\leq e$ and $\dim (C)\leq e$. Our claim follows.\fini

\begin{cor}\label{apextc}
 Let $c$ and $e$ be two integers. Assume that locally at any prime $\ip$ such that $\dim R/\ip =e+1$, $M_\ip$ is Cohen-Macaulay of codimension $c$, $N_\ip$ is Cohen-Macaulay, and  $\codim_{R_\ip}
(M_\ip\otimes_{R_\ip} N_\ip )=\codim_{R_\ip} (M_\ip )+\codim_{R_\ip} (N_\ip )$. Then for any $p\geq e+1$,
$$
\begin{array}{rcl}
a_{p}(\ext^{c}_{R}(M,N))&\leq &\max_{i}\{ a_{p+i}(\ext^{c}_{R}(M,R))+b_{i}(N)\}\\
&\leq &\reg (N)+\max_{i\geq p}\{ a_{i}(\ext^{c}_{R}(M,R))+i\} -p.\\
\end{array}
$$
\end{cor}

{\it Proof.} By Corollary \ref{reg2E}, $a_{p}(\ext^{c}_{R}(M,N))
\leq a_{p}(\ext^{c}_{R}(M,R)\otimes_R N)$ for $p\geq e+1$. As $R$ is Gorenstein, it follows that, locally at each $\ip$ such that $\dim R/\ip =e+1$,
$\ext^{c}_{R}(M,R)$ and $N$ are Cohen-Macaulay modules that intersects properly. Therefore
$\dim  (\tor_i^R(\ext^{c}_{R}(M,R),N))\leq e$, for $i>0$ (see e. g. \cite[1.8]{Ch2}) and so \cite[5.11 (i)]{Ch2} implies that
$a_{p}(\ext^{c}_{R}(M,R)\otimes_R N)\leq \max_{i}\{ a_{p+i}(\ext^{c}_{R}(M,R))+b_{i}(N)\}$ for
$p\geq e-1$. The conclusion follows.\fini

\begin{thm}\label{regextpi}

1) If $\dim (M\otimes_{R}N)\leq 1$, then
$$
\max_{j}\{ \reg (\ext^{j}_{R}(M,N))+j\} = \reg (N)-\indeg (M).
$$

2) Let $c$ be an integer. Assume that $\tor_{1}^{R}(M,N)$ is supported in dimension at most $1$  and that, for any $\ip$ such that $\dim (R/\ip )\geq 2$, $M_\ip \otimes_{R_\ip} N_\ip$ is Cohen-Macaulay and
$\codim_{R_\ip} (M_\ip )=c$.  Then,

(i) If $\reg (\ext^{c}_{R}(M,N))+c\leq \reg (N)-\indeg (M)$, then
$$
\max_{j}\{ \reg (\ext^{j}_{R}(M,N))+j\} = \reg (N)-\indeg (M).
$$

(ii) If $\reg (\ext^{c}_{R}(M,N))+c>\reg (N)-\indeg (M)$, then

$\reg (\ext^{j}_{R}(M,N))+j\leq \reg (N)-\indeg (M),\quad \forall j<c$,

$\reg (\ext^{c}_{R}(M,N))+c\leq \reg (N)+\reg (\ext^{c}_{R}(M,R))+c$, and

$\max_{j>c}\{ \reg (\ext^{j}_{R}(M,N))+j\} = \reg (\ext^{c}_{R}(M,N))+c-1$.

If further $\tor_{1}^{R}(M,N)$ is supported in dimension at most $0$, then
$$
\reg (\ext^{c}_{R}(M,N))=\reg (\ext^{c}_{R}(M,R)\otimes_{R}N)\leq\reg (\ext^{c}_{R}(M,R))+\reg (N).
$$

\end{thm}

{\it Proof.} Set $E^{j}:=\ext^{j}_{R}(M,N)$. Our hypotheses  in 1) implies that $E^j$ has dimension
at most 1 for all $j$. Recall that over a localization $R_\ip$, the following are equivalent : 

a) $\Tor_1^{R_\ip}(P,Q)=0$ and $P\otimes_{R_\ip}Q$ is Cohen-Macaulay,

b) $P$ and $Q$ are Cohen-Macaulay and $\codim_{R_\ip}P\otimes_{R_\ip}Q=\codim_{R_\ip}P+\codim_{R_\ip}Q$.

Notice that if $M$ is Cohen-Macaulay of projective dimension $p$, then $\Ext^j_R(M,N)\simeq \tor_{p-j}^R(\ext^p_R(M,R),N)$  and that $\ext^p_R(M,R)$ is Cohen-Macaulay with same support as $M$. This implies that
$E^j$ has dimension
at most 1 for all $j\not= c$ in 2) . By Lemma \ref{ssHE} the spectral sequence with $E_2$
term

$$
\xymatrix{
\cdots&H^{0}_{\im}(E^{c-1})&H^{0}_{\im}(E^{c})&H^{0}_{\im}(E^{c+1})\ar[ddl]&H^{0}_{\im}(E^{c+2})
\ar@{-->}[dddll]&H^{0}_{\im}(E^{c+3})\ar@{..>}[ddddlll]\\
 \cdots&H^{1}_{\im}(E^{c-1})&H^{1}_{\im}(E^{c})&H^{1}_{\im}(E^{c+1})\ar[ddl]&H^{1}_{\im}(E^{c+2})
 \ar@{-->}[dddll]&\cdots\\
 \cdots&0&H^{2}_{\im}(E^{c})&0&\cdots&\cdots\\
  \cdots&0&H^{3}_{\im}(E^{c})&0&\cdots&\cdots\\
  \cdots&0&H^{4}_{\im}(E^{c})&0&\cdots&\cdots\\
 }
 $$
 as well as the as the one with $E_1$ term  $H^{i}_{\im}(N)\otimes_R \hom_{R}(F_{j}^{M},R)$ abuts 
 to $H^{\bullet}_{\im}(M,N)$. (In this diagram, the dotted arrows show direction of maps at steps 3 or 4 of the spectral sequence --notice that the target of each map is
a quotient of a local cohomology module of $\ext^{c}_{R}(M,N)$).

 Set $b'_{i}(M):=\indeg (\tor_{i}^{R}(M,k))$ and $B_{j}:=\max_{i}\{ a_{i}(N)-b'_{j-i}(M)\}$. Notice
 that  $b_{0}'(M)=\indeg (M)$ and $b'_{\ell}(M)\geq b'_{0}(M)+\ell$ and therefore
 $B_{j}\leq \max_{i}\{ a_{i}(N)-\indeg (M)-j+i\} \leq \reg (N)-\indeg (M)-j$. Comparing the two spectral sequences, for any $j\leq c$, one has 
$
a_{0}(\ext^{j}_{R}(M,N))\leq B_{j}
$
and
$
a_{1}(\ext^{j}_{R}(M,N))\leq B_{j+1}
$.
Therefore, for $j\leq c$,
$$
\{ \; a_{0}(\ext^{j}_{R}(M,N))+j\leq \reg (N)-\indeg (M)\quad \& \quad
a_{1}(\ext^{j}_{R}(M,N))+1+j\leq \reg (N)-\indeg (M)\; \} \quad (*).
$$
This shows that, in both part 1 and part 2, $\reg (\ext^{j}_{R}(M,N))+j\leq \reg (N)-\indeg (M),\ \forall j<c$.

For $j>c$, again by comparing the two spectral sequences,
$$
\begin{array}{rcl}
a_{0}(\ext^{j}_{R}(M,N))+j&\leq& \max\{ B_{j}+j,a_{j-c+1}(\ext^{c}_{R}(M,N))+j\}\\
&\leq& \max\{ \reg (N)-\indeg (M),\reg (\ext^{c}_{R}(M,N))+c-1\}\\
\end{array}
$$
$$
\begin{array}{rcl}
a_{1}(\ext^{j}_{R}(M,N))+1+j&\leq& \max\{ B_{j+1}+j+1,a_{j-c+2}(\ext^{c}_{R}(M,N))+j+1\}\\
&\leq&\max\{  \reg (N)-\indeg (M),\reg (\ext^{c}_{R}(M,N))+c-1\}\\
\end{array}$$

The above estimates and $(*)$ imply that in part 1, $\max_j \{ \reg (\ext^{j}_{R}(M,N))+j \} \leq \reg (N)-\indeg (M)$. In part 2, the two above estimates show that 
$$
\max_{j>c} \{ \reg (\ext^{j}_{R}(M,N))+j \} \leq
\max\{  \reg (N)-\indeg (M),\reg (\ext^{c}_{R}(M,N))+c-1\} .
$$ 


Furthermore, recall that by Theorem \ref{greg} 5), there exists $p$ such that 
$\fin (H^{p}_{\im}(M,N))=\reg (N)-\indeg (M)-p$. 

In part 1), this shows that $a_{0}(\ext^{p}_{R}(M,N))=\reg (N)-\indeg (M)-p$ or
$a_{1}(\ext^{p-1}_{R}(M,N))=\reg (N)-\indeg (M)-p$. It implies that either $\reg (\ext^{p}_{R}(M,N))\geq \reg (N)-\indeg (M)-p$ or $\reg (\ext^{p-1}_{R}(M,N))\geq \reg (N)-\indeg (M)-p+1$. This finishes the proof of 1).

In part 2) (i), this shows the same inequalities if $p\leq c$. If $p>c$, it shows that $a_{0}(\ext^{p}_{R}(M,N))\geq \reg (N)-\indeg (M)-p$ or $a_{1}(\ext^{p-1}_{R}(M,N))\geq \reg (N)-\indeg (M)-p$ or
$a_{p-c}(\ext^{c}_{R}(M,N))\geq \reg (N)-\indeg (M)-p$. This proves 2) (i).

Finally, in part 2) (ii), choose $i$ such that $\reg (\ext^{c}_{R}(M,N))=a_{i}(\ext^{c}_{R}(M,N))+i$. 
Notice that $i\geq 2$ by $(*)$. Hence it follows from  Corollary \ref{apextc} that $\reg (\ext^{c}_{R}(M,N))\leq \reg (\ext^{c}_{R}(M,R))+\reg (N)$. 
Let $\mu :=a_{i}(\ext^{c}_{R}(M,N))$. As $\mu >\reg (N)-\indeg (M)-i-c$ by hypothesis, $H^{1}_{\im}(\ext^{c+i-2}_{R}(M,N))_\mu=H^{0}_{\im}(\ext^{c+i-1}_{R}(M,N))_\mu =0$ implies that $H^{i}_{\im}(\ext^{c}_{R}(M,N))_\mu = 0$. Hence, either $a_{1}(\ext^{c+i-2}_{R}(M,N))\geq 
\mu$ or $a_{0}(\ext^{c+i-1}_{R}(M,N))\geq 
\mu$. This shows that $\max_{j>c}\{ \reg (\ext^{j}_{R}(M,N))+j\} = \reg (\ext^{c}_{R}(M,N))+c-1$.

In view of $(*)$, Corollary \ref{reg2E} implies the last statement of part 2) (ii). \fini

 \bco
 Under the hypotheses of Theorem \ref{regextpi} 1) or 2) (i), 
 $$
 \begin{array}{rl}
 \max_{j}\{ \reg (\ext^{j}_{R}(M,N))+j\} &\!\!\!\!  -\min_{j}\{ \indeg (\ext^{j}_{R}(M,N))+j\}\\
 &= \reg (M)-\indeg (M)+\reg (N)-\indeg (N)\\
 &= \reg_R (M,M)+\reg_R (N,N)= \reg_R (M,N)+\reg_R (N,M).\\
 \end{array}
 $$
 \eco
 
 We collect below some results about the regularity of the Ext module playing a key 
 role in the estimates of Theorem \ref{regextpi}, in the case $M$ is cyclic.
 
 \bt
 Let $I$ be a homogeneous proper ideal of $R$, $c:=\codim_{R} (I)$, $J$ the unmixed part of $I$ and $X:=\proj (R/J)$. Then,
 $$
 \reg (\ext^{c}_{R}(R/I,R))+c=0
 $$
 if either
 
 (i) $X$ is empty or arithmetically Cohen-Macaulay,
 
 (ii) $\dim (X)=1$ and $X$ is reduced,
 
 (iii) $X$ is smooth, $\chara (R_0)\geq \dim (X)$ and $X$ lifts to $W_2(R_0)$,
 
 (iv) $\chara (R_0)=0$ and $X$ has at most isolated irrational singularities.
 \et
 
 {\it Proof.} If $R/J$ is Cohen-Macaulay, then $\ext^{c}_{R}(R/J,R)$ has the dual of a minimal free $R$-resolution of $R/J$ as minimal free $R$-resolution, from which (i) directly follows. 
 
Notice that $\reg (\ext^{c}_{R}(R/I,R))+c=0$ is equivalent to $\reg (\om_{R/J})=\dim (R/J)$, 
 with $\om_{R/J}:=\ext^{c}_{R}(R/J,\om_{R})$ and that $\reg (\om_{R/J})\geq \dim (R/J)$ (see \cite[4.2]{CU}). Part (ii) is proved for instance in the first line of the proof of \cite[1.7]{CU}, (iii) follows from \cite[2.8]{DI} and (iv) follows from \cite[1.3]{CU}.\fini
  
 \brm
 In the case $\dim (X)=1$ above, the hypothesis that $X$ is reduced or ACM is essential, as 
 the following example shows (see \cite[13.5]{Ch1}). Let $R:=k[x,y,z,t]$ and $I:=(x^n t-y^n z)+(z,t)^n$ for
 $n\geq 2$. Then $I$ is unmixed, $\proj (R/I)$ is locally a complete intersection and $\reg (\ext^{2}_{R}(R/I,R))+2=(n-1)^2$.
 \erm

\end{document}